\documentclass[%
nofootinbib,
 amsmath,amssymb,notitlepage
]{revtex4-1}
\pdfoutput=1

\usepackage{graphicx}
\usepackage{dcolumn}
\usepackage{bm}

\usepackage{physics}
\usepackage{wrapfig}
\usepackage{hanging}
\usepackage{hhline}

\DeclareMathOperator{\Mod}{mod}

\begin{document}

\widowpenalty10000
\clubpenalty10000

\title{A 5-Dimensional \textit{Tonnetz} for Nearly Symmetric Hexachords}
\thanks{I thank Professor Suzannah Clark for discussions during the preparation of this paper.}%

\author{Vaibhav Mohanty}
\email{E-mail: mohanty@college.harvard.edu}
\affiliation{%
 Quincy House, Harvard University, Cambridge, MA 02138
}%

\count\footins = 10000


\begin{abstract}
The standard 2-dimensional \textit{Tonnetz} describes parsimonious voice-leading connections between major and minor triads as the 3-dimensional \textit{Tonnetz} does for dominant seventh and half-diminished seventh chords. In this paper, I present a geometric model for a 5-dimensional \textit{Tonnetz} for parsimonious voice-leading between nearly symmetric hexachords of the mystic-Wozzeck genus. Cartesian coordinates for points on this discretized grid, generalized coordinate collections for 5-simplices corresponding to mystic and Wozzeck chords, and the geometric nearest-neighbors of a selected chord are derived.
\end{abstract}

\maketitle

\section{Introduction}
In this paper, I constrct a 5-dimensional \textit{Tonnetze} for nearly symmetric hexachords, known as the mystic and Wozzeck chords. Cohn (1996) describes that it should be possible for chords of the $T_n/T_nI$ set class of 6-34 to exhibit voice-leading parsimony in ways similar to the major and minor triads of 3-11. He showed that the major and minor triads, which can be viewed as perturbations of the (symmetric) augmented triad, exhibit smooth chord transitions which are achieved via two separate voice-leading regions: the hexatonic region and the Weitzmann waterbug region (Cohn 2012). Separate Neo-Riemannian transformations exist within each of these regions, and all of the transformations contained in the union of these two sets, with the exception of the hexatonic pole relation $\vb{H}$, are visualizable on the well-known 2D \textit{Tonnetz}.   

The Boretz spider region and Childs' (1998) octatonic region similarly contain Neo-Riemannian transformations that appropriately describe parsimonious voice-leading between dominant seventh and half-diminished seventh chords, which are perturbations of the symmetric fully diminished chord. These chords are visually reprented as regular tetrahedra (3-simplices) in Gollin's (1998) 3D \textit{Tonnetz}, and the Neo-Riemannian transformations appear naturally as nearest-neighbor relations.

The author has previously developed the dodecatonic and centipede regions for voice-leading between chords obtained from the perturbation of the symmetric whole-tone scale: the mystic and Wozzeck chords (Mohanty 2018). In this paper, I mathematically construct the 5D \textit{Tonnetz} for voice-leading between mystic and Wozzeck chords. First, the positions of the pitch classes in $\mathbb{R}^5$ must be established. From there, the Wozzeck and mystic chords can be defined by the positions of their vertices. Lastly, nearest-neighbor chords can be found, and Neo-Riemannian transformations from the previous work (Mohanty 2018) can be assigned where applicable. 

\section{The coordinate space}
In this section, I construct the geometric positions and identities of individual pitch classes in 5D space. Like the 2D and 3D \textit{Tonnetze}, equally spaced points in the space represent pitch classes, and simplices bounded by $n-1$ vertices in $\mathbb{R}^n$ correspond to nearly symmetric chords.
\subsection{The {\normalfont \textbf{Tonnetz}} basis of $\mathbb{R}^5$}
In both the 2D and 3D note spaces, the axes along which individual pitch classes lie are not mutually orthogonal. Writing the unit vectors pointing along these axes as $\{\vu{q}_1, \dots, \vu{q}_{N-1}\}$, where $N$ is the cardinality of the chord, it is easy to see that the $i$-th and $j$-th unit vectors will satisfy
\begin{equation}
\vu{q}_i \vdot \vu{q}_j = \frac{1}{2}.
\label{relation}
\end{equation}
I generalize this relation to the $N = 6$ case so that all 5 axes in the 5-dimensional note space will be oriented at $60$ degrees with respect to one another. Imposing the conditions in eq. (\ref{relation}) separately on the 5 axes, we find that the unit vectors $\{\vu{q}_i\}$ can be written in the Cartesian basis $\{\vu{e}_i\}$ as
\begin{equation}
\begin{split}
\vu{q}_1 = \begin{pmatrix}
1 \\
0 \\
0 \\
0 \\
0
\end{pmatrix}, \quad
\vu{q}_2 = \begin{pmatrix}
1/2 \\
\sqrt{3}/2 \\
0 \\
0 \\
0
\end{pmatrix}, \quad
\vu{q}_3 = \begin{pmatrix}
1/2 \\
1/(2\sqrt{3}) \\
\sqrt{2/3} \\
0 \\
0
\end{pmatrix}, \\\quad
\vu{q}_4 = \begin{pmatrix}
1/2 \\
1/(2\sqrt{3}) \\
1/(2\sqrt{6}) \\
\sqrt{5/8} \\
0
\end{pmatrix}, \quad
\vu{q}_5 = \begin{pmatrix}
1/2 \\
1/(2\sqrt{3}) \\
1/(2\sqrt{6}) \\
1/(2\sqrt{10}) \\
\sqrt{3/5}
\end{pmatrix}. \quad
\end{split}
\end{equation}
Any vector $[\vb{v}]_{\vb{q}} \in \mathbb{R}^{5}$ in the \textit{Tonnetz} basis $\{\vu{q}_i\}$ can be represented in the Cartesian basis as $\vb{v} = U[\vb{v}]_{\vb{q}}$ by the unitary transformation
\begin{equation}
U = \begin{pmatrix}
\vu{q}_1 & \vu{q}_2 & \vu{q}_3 & \vu{q}_4 & \vu{q}_5
\end{pmatrix} = \begin{pmatrix}
1 & 1/2 & 1/2 & 1/2 & 1/2 \\
0 & \sqrt{3}/2 & 1/(2/\sqrt{3}) & 1/(2/\sqrt{3}) & 1/(2/\sqrt{3}) \\
0 & 0 & \sqrt{2/3} & 1/(2\sqrt{6}) & 1/(2\sqrt{6}) \\
0 & 0 & 0 & \sqrt{5/8} & 1/(2\sqrt{10}) \\
0 & 0 & 0 & 0 & \sqrt{3/5}
\end{pmatrix}.
\end{equation}
\subsection{Pitches in the coordinate space}
Let $S$ denote the set of \textbf{tones} in the \textit{Tonnetz} coordinate space; in particular, $S$ includes all linear combinations of the \textit{Tonnetz} basis vectors $\{\vb{q}_i\}$ with integer coefficients. That is,
\begin{equation}
S = \{i\vu{q}_1 + j\vu{q}_2 + k\vu{q}_3 + \ell\vu{q}_4 + m\vu{q}_5 \:|\: i,j,k,\ell,m \in \mathbb{Z}\}.
\end{equation}
We define a map $\varphi : S \rightarrow \mathbb{Z}_{12}$ such that $\varphi(\vb{s})$ for $\vb{s} \in S$ returns an integer $\varphi(\vb{s}) \in \{0,\dots,11\}$ that corresponds to a particular pitch class $\{C,\dots,B\}$, and assignment is inherently arbitrary. However, throughout this paper, I use the standard convention of $0 = C$, $1 = C\sharp$, etc. I will also use an ordered pair of integers $(i,j,k,\ell,m)$ to represent elements of $S$ instead of the standard column vector; this notation should not be confused with my notation for row vectors, which are not written with commas.

Now, I explicitly construct $\varphi$ by following the conventions of the 2D and 3D \textit{Tonnetze}. We can succinctly state that
\begin{equation}
\varphi(i\vu{q}_1 + j\vu{q}_2 + k\vu{q}_3 + \ell\vu{q}_4 + m\vu{q}_5) = \Mod_{12}(4i + 8j + 10k + \ell + 6m)
\end{equation}
where $\Mod_{12} : \mathbb{Z} \rightarrow \mathbb{Z}_{12}$ returns the remainder of the argument divided by 12. From this definition, one may see that, starting at the origin, the notes along the $\vu{q}_1$ in the positiive direction are $C$, $E$, $G\sharp$, etc. The notes along the $\vu{q}_2$ in the positive direction are $C$, $A\flat$, $E$, etc. Similar logic can be applied to the other three axes. Notes that are not along any axis are determined simply by linearity.

\section{Nearly symmetric hexachords}
Now that the coordinate space has been constructed precisely, I now introduce the geometric definitions of the mystic and Wozzeck chords. As described in section 1.3, the mystic and Wozzeck chords are inversionally related nearly symmetric hexachords, and I will show that the particular definition of $\varphi$ in the previous section has been provided so that each mystic chord and each Wozzeck chord forms a 5-simplex in $\mathbb{R}^5$.
\subsection{Wozzeck chords}
A Wozzeck chord is obtained from the downward perturbation of any tone in a whole-tone scale and will be denoted with a $(+)$ symbol such that ``$C$ Wozzeck'' can be written as $C+$. By the convention presented in an earlier work (Mohanty 2018), a Wozzeck chord will be labeled by the lower of the two tones comprising a minor 2nd. This is to say that $C+$ is the collection of pitch classes $\{ C, D\flat, E, F\sharp, G\sharp, B\flat\}$. 

In the coordinate space defined in the previous section, a Wozzeck chord which has its root at the point $(i,j,k,\ell,m)$ is given by the collection of tones $\{(i,j,k,\ell,m),(i+1,j,k,\ell,m),(i,j+1,k,\ell,m),(i,j,k+1,\ell,m),(i,j,k,\ell+1,m),(i,j,k,\ell,m+1)\}$. This corresponds to the collection of vertices of a 5-simplex in $\mathbb{R}^5$ with orientation we will refer to as $(+)$.
\subsection{Mystic chords}
A mystic chord is given by a upward perturbation of a tone within the whole-tone scale; these chords are denoted with the $(-)$ symbol. Thus, $C-$ refers to the ``$C$ mystic'' chord and is comprised of the pitch classes $\{C, D\flat, E\flat, F, G, A\}$. A mystic chord is given by the collection of 6 tones $\{(i,j,k,\ell,m), (i + 1, j, k, \ell, m), (i + 1, j - 1, k, \ell, m), (i + 1, j, 
k - 1, \ell, m), (i + 1, j, k, \ell - 1, m), (i + 1, j, k, \ell, m - 1)\}$. These points directly correspond to the set of vertices of a 5-simplex with orientation opposite to that of the Wozzeck chords---the $(-)$ orientation.
\subsection{Duality in the mystic-Wozzeck genus }
As major and minor triads---represented by triangles (or 2-simplices) in the 2D \textit{Tonnetz}---have opposite graphical orientations, the dominant seventh and half-diminished seventh chords in Gollin's (1998) 3D \textit{Tonnetz} also are ``upside down'' images of each other. This notion of orientation is well-described mathematically and can easily be obtained by comparing the set of $\mathbb{R}^5$ coordinates with the general forms of the Wozzeck collection from section 3.1 and the mystic collection from section 3.2. A nearly symmetric hexachord can only be represented in one of the two inversionally related forms, so the notions of $(+)$ or $(-)$ orientation holds for the mystic-Wozzeck genus as it does for the major and minor triads as well as the Tristan genus. The mathematical notion of orientation, which is a signed quantity, preserves this analogy as well.
\section{Neighbors in the 5D \textit{Tonnetz}}
As described by Cohn (2012), the nearly symmetric hexachords exhibit parsimonious voice leading, and small-displacement chord transitions are fully described by a set of Neo-Riemannian transformations, which are defined in the author's previous work (Mohanty 2018). The table of these Neo-Riemannian transformations and the result of applyinng these transformations to $C+$ are displayed in \textbf{Table 1}.
\begin{table}[t]
\small
\caption{Summary of Neo-Riemannian transformations for nearly symmetric hexachords and results of operations on $C+$.}
\centering
\begin{tabular}{|c|c|c|}
\hline
Starting Chord & Transformation & Resulting Chord\\
\hhline{|=|=|=|}
\multicolumn{3}{|c|}{For $P_{0,1}$-related chords}\\ \hline
$C+$ & $\vb{R}^{**}$ & $D\sharp-$  \\\hline
\multicolumn{3}{|c|}{For $P_{2,0}$-related chords}\\ \hline
$C+$ & $\vb{S^{A(3)}}$ & $B-$  \\\hline
$C+$ &  $\vb{S^{A(5)}}$ & $G-$  \\\hline
$C+$ & $\vb{S^{F}}$ & $A-$  \\\hline
$C+$ & $\vb{S^{W(1)}}$ & $C\sharp-$  \\\hline
$C+$ & $\vb{S^{W(3)}}$ & $F-$  \\\hline
\multicolumn{3}{|c|}{For $P_{n-2,0}$-related chords}\\ \hline
$C+$ & $\vb{S^1}$ & $C-$  \\\hline
$C+$ & $\vb{S^{3(A)}}$ & $A\sharp-$  \\\hline
$C+$ & $\vb{S^{3(W)}}$ & $E-$  \\\hline
$C+$ & $\vb{S^{5(A)}}$ & $F\sharp-$  \\\hline
$C+$ & $\vb{S^{5(F)}}$ & $G\sharp-$  \\\hline
\multicolumn{3}{|c|}{Polar Relation}\\ \hline
$C+$ & $\vb{Z}$ & $D-$  \\\hline
\end{tabular}
\end{table}
Starting with some arbitary chord on the 2D \textit{Tonnetz}, it is easy to see that applying the standard triadic Neo-Riemannian transformations $\vb{R}$, $\vb{P}$, $\vb{L}$, $\vb{S}$, and $\vb{N}$ to the starting chord result in a chord that shares either an edge or a corner with the starting chord. In the 2D and 3D \textit{Tonnetze}, not all of the neighboring corner or edge chords are represented by the above transformations, but for the 5D \textit{Tonnetz} every neighbor has an associated transformation. The polar relation $\vb{H}$ does not---and should---not share any common tones with the starting chord, so it is not a neighbor. Examining Gollin's (1998) diagrams, it is clear that the same rule holds for the 3D \textit{Tonnetz}; all of the well-defined Neo-Riemannian transformations except the octatonic pole $\vb{O}$ transformation correspond either to an edge-preserving or corner-preserving neighbor chord. One may expect the rule to hold for the 5D \textit{Tonnetz}, and indeed it does, as I will show.

Since the full 5-dimensional space cannot be directly visualized in spatial coordinates, I have produced several reduced images in \textbf{Figures 1} through \textbf{7}. In each figure 1-7, the central chord, which appears as a hexagon with several diagonal lines, is the arbitrarily chosen $C+$. This hexagon represents an orthographic projection of the $C+$ 5-simplex onto 2 dimensions, and every solid line represents the edge of a 5-simplex in the $\mathbb{R}^5$ \textit{Tonnetz} space. Despite varying lengths of the solid lines in this projection, each line represents the same $\mathbb{R}^5$ distance, which is precisely unit distance using the standard Euclidean metric.

The permutation of the vertex labels of a given chord in different \textbf{Figures 1} through \textbf{7} allow for easy visualization of neighbors. A 5-simplex has 15 edges and 6 corners, so ``true'' picture of the 5D \textit{Tonnetz} is a simultaneous superposition of all 7 panels shown in \textbf{Figures 1} through \textbf{7}. In the figure, the Neo-Riemannian transformation relating $C+$ and the neighboring chord---if such a transformation is well-defined---is given in bold next to the neighboring chord.

The rules for the chord neighbors shown in \textbf{Figures 1} through \textbf{7} generally hold for any Wozzeck chord, and the neighbors for any mystic chord can be quickly deduced by symmetry properties. A limitation of the 5D \textit{Tonnetz} is the inevitable fact that the entire \textit{Tonnetz} cannot be visualized with accurate representation of all spatial dimensions simultaneously. The orthographic projections used in this paper, moreover, cause the Wozzeck and mystic chords to appear geometrically as identical objects, whereas the 2D and 3D \textit{Tonnetz} clearly distinguish between chords of opposite quality by clearly displaying orientation of simplices. For the 5D \textit{Tonnetz}, the reader must actively examine the identities of the vertex pitch classes of a particular chord to parse whether the examined chord is a mystic or Wozzeck chord.
\section{Conclusion}
In this article, I have presented an explicit construction of the 5D \textit{Tonnetz} for voice-leading between nearly symmetric hexachords. As discussed in previous work (Mohanty 2018; Cohn 2012), Mystic and Wozzeck chords obey voice-leading rules similar to those for the major-minor triadic complex as well as the Tristan genus. The 5D \textit{Tonnetz} presented here is intended as an analogy and extension of the 2D and 3D \textit{Tonnetze} to the remaining class of perturbatively constructed chords of cardinality $n = 6$. As Cohn's (1996) hexatonic and Weitzmann waterbug regions define Neo-Riemannian transformations for major and minor chords that can be represented on the 2D \textit{Tonnetz}, Childs's (1998) octatonic and Boretz spider regions present Neo-Riemannian transformations that are used to voice-lead between dominant seventh and half-diminished seventh chords that can be represented on Gollin's (1998) 3D \textit{Tonnetz}. For mystic and Wozzeck chords, the dodecatonic and centipede regions (Mohanty 2018) are comprised of Neo-Riemannian transformations that can be depicted within the 5D \textit{Tonnetz} presented here.
\section*{References}
{\small
\begin{hangparas}{1cm}{1}
\hangindent=1cm 
Childs, Adrian P. 1998. ``Moving beyond Neo-Riemannian Triads: Exploring a Transformational Model for Seventh Chords.'' \textit{Journal of Music Theory} 42, no. 2: 181-193.

Cohn, Richard. 1996. ``Maximally Smooth Cycles, Hexatonic Systems, and the Analysis of Late-Romantic Triadic Progressions.'' \textit{Music Analysis} 15, no. 1: 9-40.

Cohn, Richard. 2012. \textit{Audacious Euphony: Chromatic Harmony and the Triad's Second Nature}. 2nd Edition. New York: Oxford University Press.

Gollin, Edward. 1998. ``Some Aspects of Three-Dimensional `\textit{Tonnetze}'.'' \textit{Journal of Music Theory} 42, no. 2: 195-206.

Mohanty, Vaibhav. 2018. ``Dodecatonic Cycles and Parsimonious Voice-Leading in the Mystic-Wozzeck Genus.'' Submitted for publication. Preprint: arXiv:1805.11087.
\end{hangparas}
}
\vfill
\section*{Appendix: Figures}
\begin{figure}[h]
  \centering
  \includegraphics[width=\textwidth]{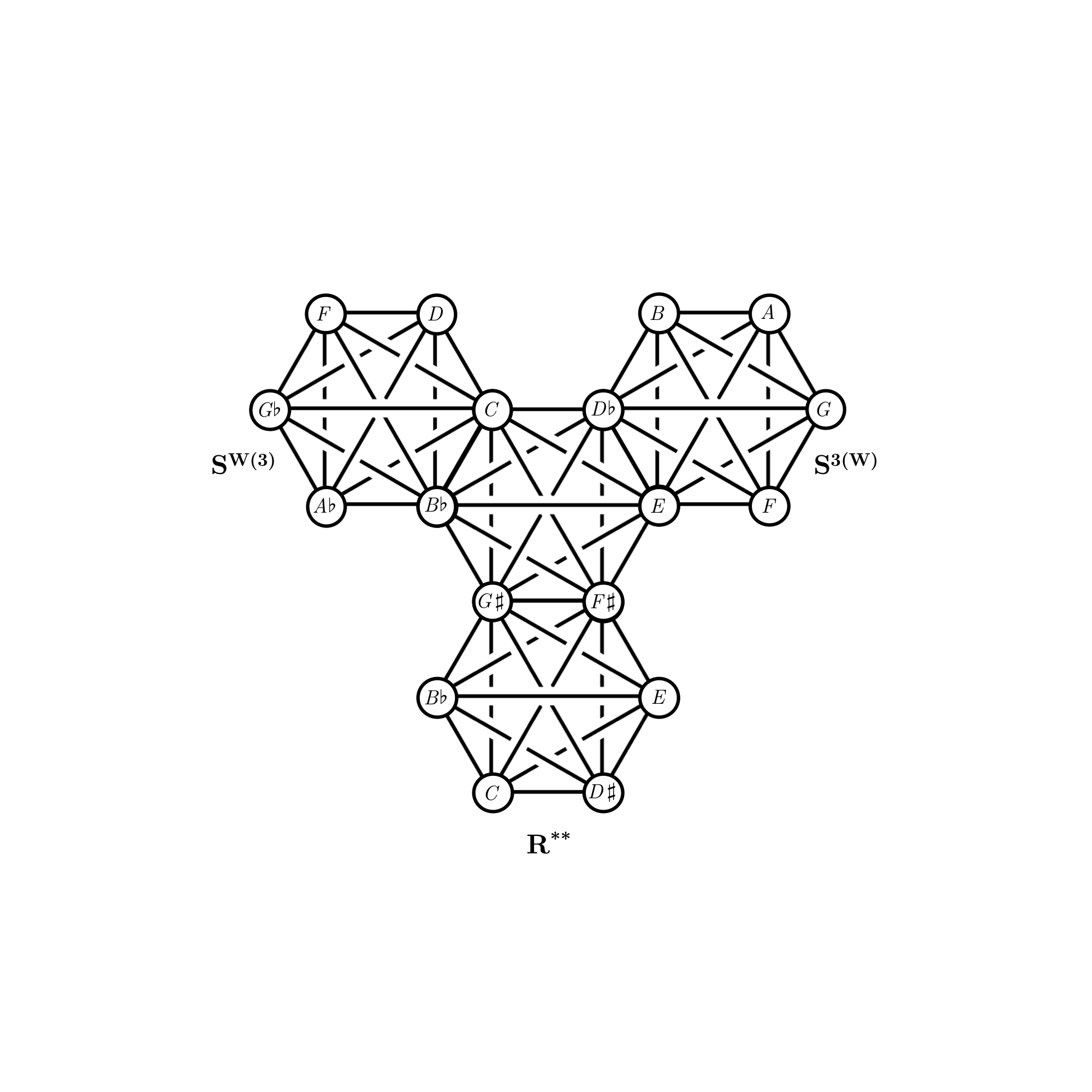}
  \caption{Edge-sharing chords in the 5D \textit{Tonnetz}. The central chord is $C+$.}
\end{figure}
\begin{figure}[h]
  \centering
  \includegraphics[width=\textwidth]{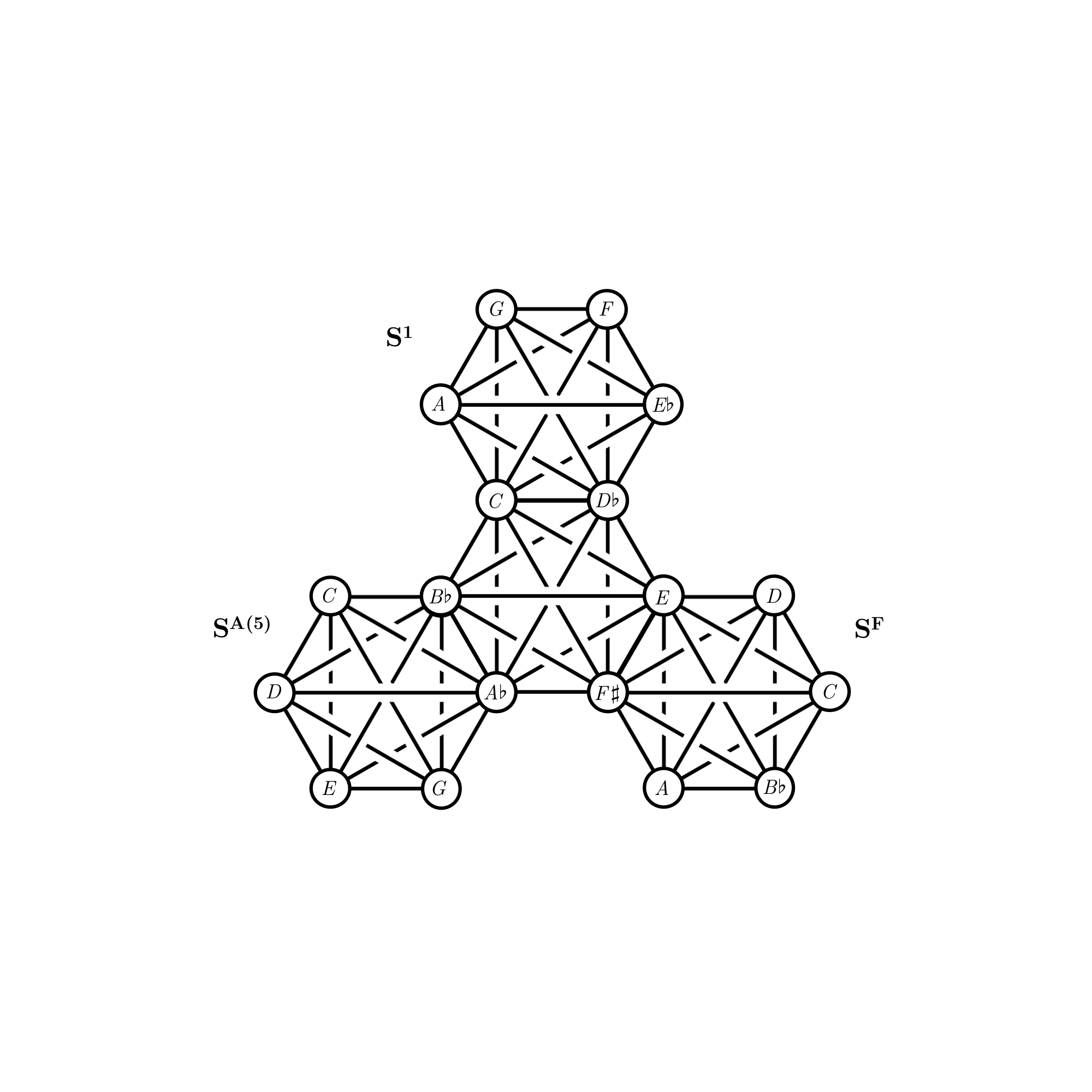}
  \caption{Edge-sharing chords in the 5D \textit{Tonnetz}. The central chord is $C+$.}
\end{figure}
\begin{figure}[h]
  \centering
  \includegraphics[width=\textwidth]{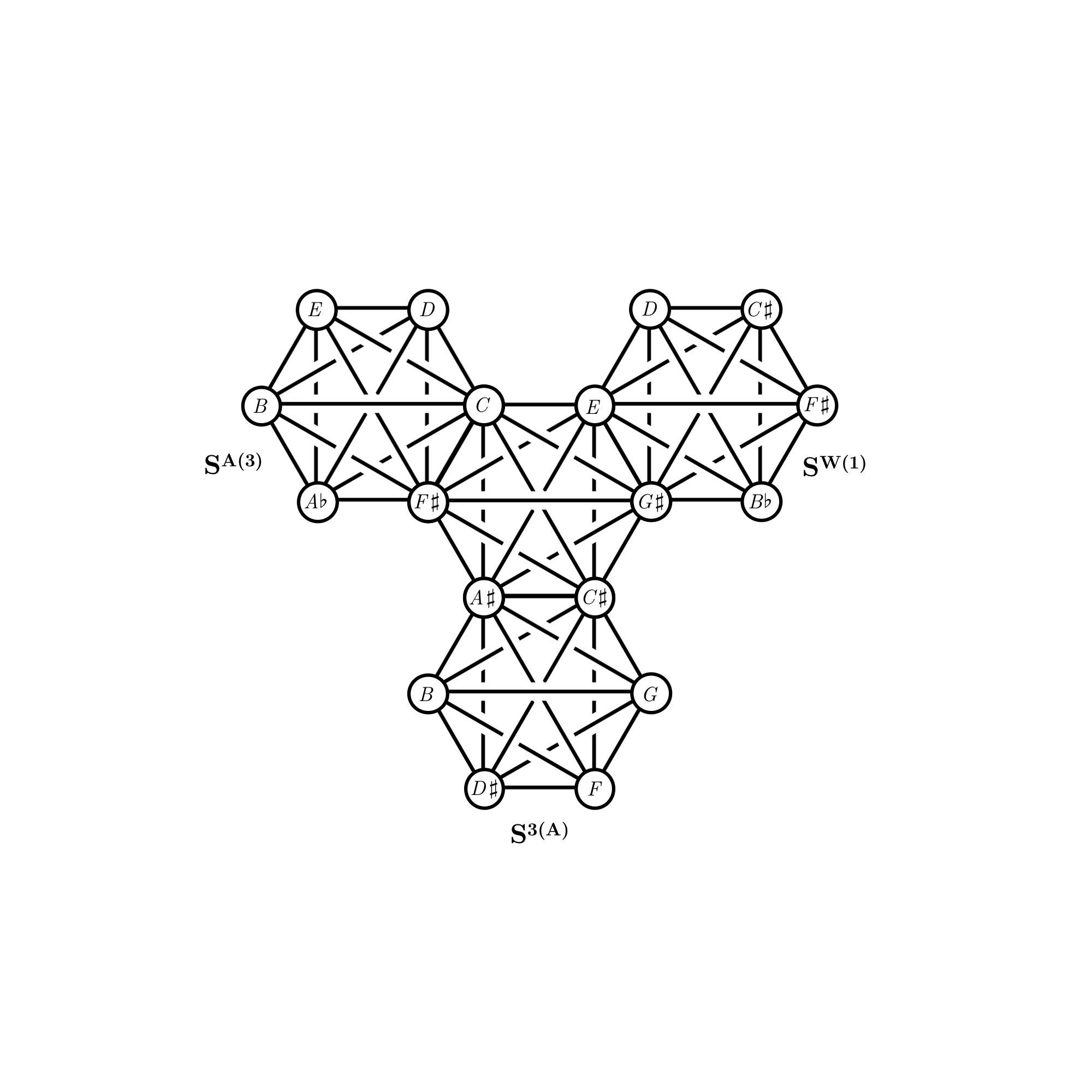}
  \caption{Edge-sharing chords in the 5D \textit{Tonnetz}. The central chord is $C+$.}
\end{figure}
\begin{figure}[h]
  \centering
  \includegraphics[width=\textwidth]{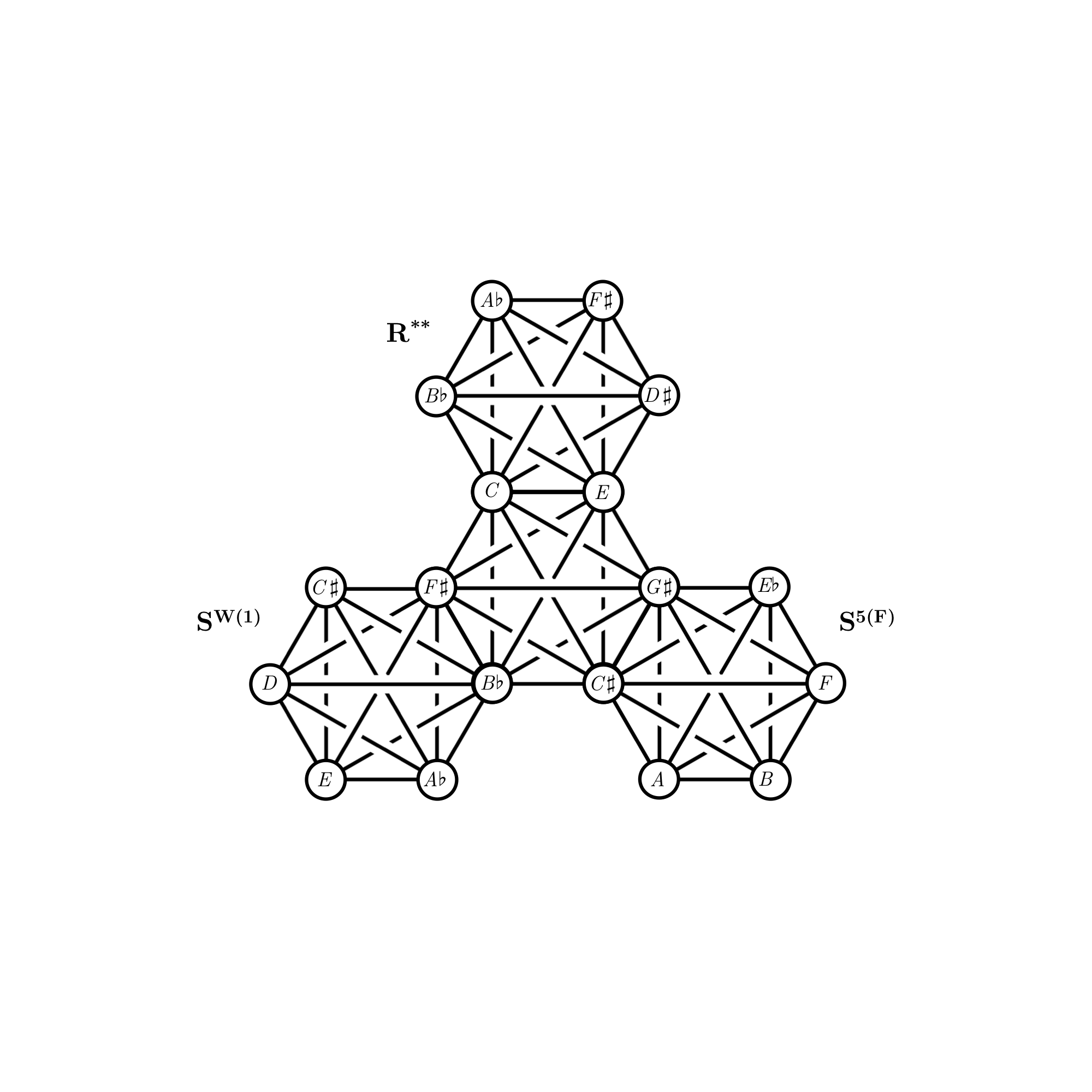}
  \caption{Edge-sharing chords in the 5D \textit{Tonnetz}. The central chord is $C+$.}
\end{figure}
\begin{figure}[h]
  \centering
  \includegraphics[width=\textwidth]{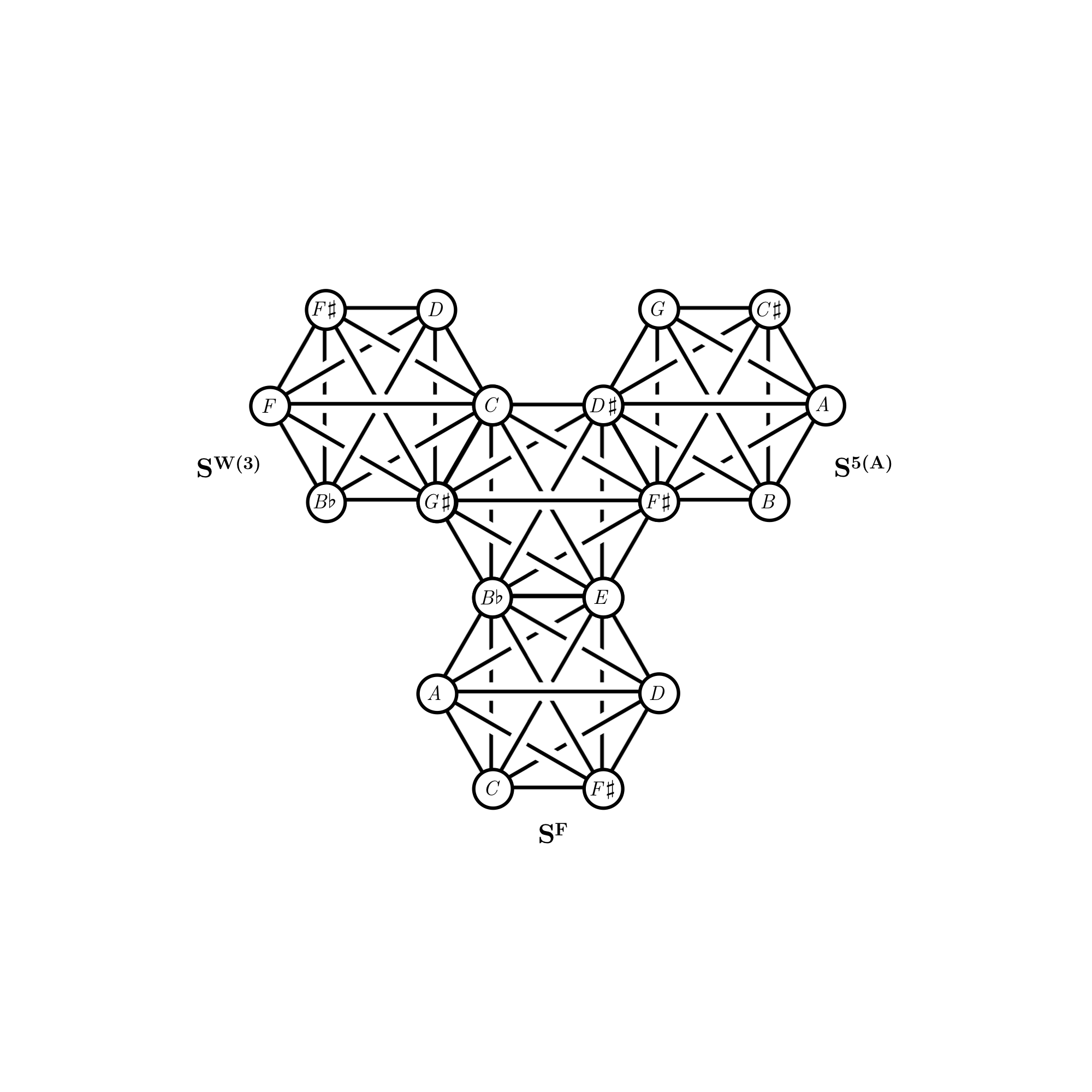}
  \caption{Edge-sharing chords in the 5D \textit{Tonnetz}. The central chord is $C+$.}
\end{figure}
\begin{figure}[h]
  \centering
  \includegraphics[width=\textwidth]{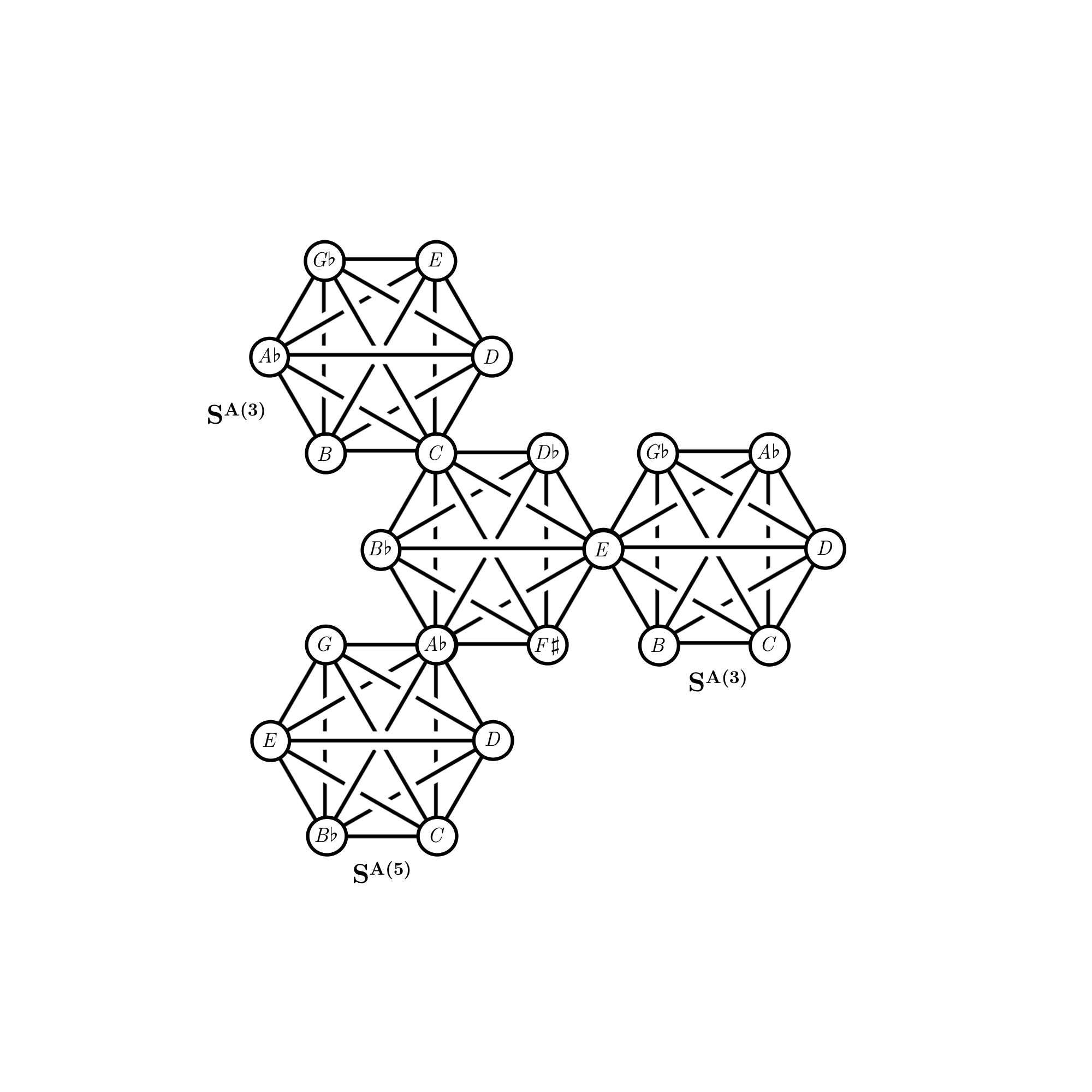}
  \caption{Corner-sharing chords in the 5D \textit{Tonnetz}. The central chord is $C+$.}
\end{figure}
\begin{figure}[h]
  \centering
  \includegraphics[width=\textwidth]{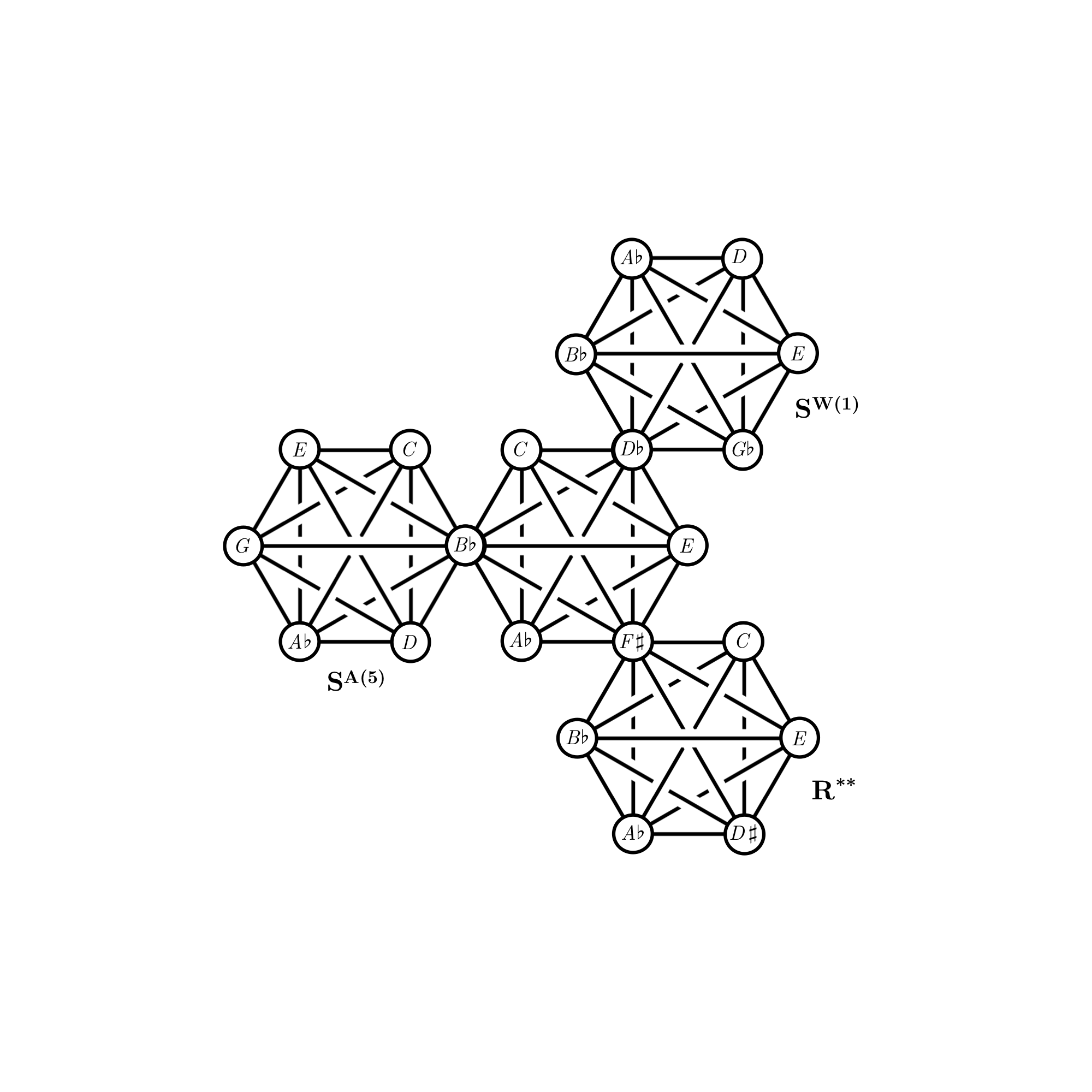}
  \caption{Corner-sharing chords in the 5D \textit{Tonnetz}. The central chord is $C+$.}
\end{figure}

\end{document}